\documentclass[12pt]{article}
\usepackage[latin2]{inputenc}
\usepackage{amsfonts}
\usepackage{amsmath}
\usepackage{enumerate}
\usepackage{amssymb}
\usepackage{theorem}
\usepackage{authblk}
\usepackage{hyperref}

\newcommand{\R}{{\ensuremath{\mathbb{R}}}}

\renewcommand{\P}{\ensuremath{\mathbb{P}}}
\renewcommand{\dj}{d\kern-0.4em\char"16\kern-0.1em}
\newcommand{\E}{\ensuremath{\mathcal{E}}}
\newcommand{\B}{\ensuremath{\mathcal{B}}}

\newcommand{\proof}{\noindent\textbf{Proof.}\ }

\newcommand{\qed}{\hfill\ensuremath{\Box}\\}

\newcommand{\HH}{{\bf (H) }}

\newcommand{\EE}{\mathbb{E}}
\newcommand{\F}{\mathcal{F}}

\newtheorem{Thm}{Theorem}[section]

\newtheorem{Lem}[Thm]{Lemma}

\newtheorem{Def}[Thm]{Definition}

{\theorembodyfont{\upshape} }
\numberwithin{equation}{section}

\title{Boundary Harnack principle for the absolute value of a one-dimensional subordinate Brownian motion killed at 0}
\author{Vanja Wagner\footnote{This work was supported by the Croatian Science Foundation under the project 3526}}
\affil{Department of Mathematics, University of Zagreb}
\date{\vspace{-5ex}}

\begin{document}
\maketitle
\begin{abstract}We prove the Harnack inequality and boundary Harnack principle for the absolute value of a one-dimensional recurrent subordinate Brownian motion killed upon hitting 0, when 0 is regular for itself and the Laplace exponent of the subordinator satisfies certain global scaling conditions. Using the conditional gauge theorem for symmetric Hunt processes we prove that the Green function of this process killed outside of some interval $(a,b)$ is comparable to the Green function of the corresponding killed subordinate Brownian motion. We also consider several properties of the compensated resolvent kernel $h$, which is harmonic for our process on $(0,\infty)$.
\end{abstract}


\section{Introduction}
Let $X=(X_t)_{t\geq 0}$ be a recurrent subordinate Brownian motion on $\mathbb R$ such that 0 is regular for itself with $\phi:(0,\infty)\to(0,\infty)$ the Laplace exponent of the corresponding subordinator. We assume that $\phi$ is a complete Bernstein function satisfying a certain global scaling condition {\bf (H)}.

The goal of this paper is to establish the Harnack inequality and boundary Harnack principle for nonnegative harmonic functions of the absolute value of process $X$ killed upon hitting $\{0\}$, denoted by $Z=(Z_t)_{t\geq0}$. In order to do so, we show that the Green function for the killed process $Z^{(a,b)}$ on a finite interval $(a,b)$, $a>0$, is comparable to the Green function of $X^{(a,b)}$. We introduce a third process $Y=(Y_t)_{t\geq 0}$ on $(0,\infty)$, obtained from $X^{(0,\infty)}$ by creation through the Feynman-Kac transform with rate equal to the killing density $\kappa_{(0,\infty)}$, see \eqref{eq:resurrect}. Process $Y$ is called the resurrected (censored) process on $(0,\infty)$ corresponding to $X$. Using the conditional gauge theorems from \cite{gauge} and the sharp two-sided Green function estimates \eqref{killed:G} for $X^{(a,b)}$ obtained in \cite{plms14} we first show that the Green functions of processes $X^{(a,b)}$ and $Y^{(a,b)}$ are comparable. In the second step we relate $Y^{(a,b)}$ and $Z^{(a,b)}$ through a Feynman-Kac transform by a discontinuous additive functional and apply the corresponding conditional gauge theorem. 

We examine a function $h:\mathbb R\to(0,\infty)$ defined by
\begin{equation*}
h(x)=\lim\limits_{q\downarrow0}(u^q(0)-u^q(x)),
\end{equation*}
which is harmonic for $Z$ on $(0,\infty)$. Here $u^q$ is the $q$-potential density of process $X$, so $h$ is sometimes called the compensated resolvent kernel. This function is often considered in relation to the local time of L\' evy processes and its properties were extensively studied in \cite{spl2006}, \cite{yano} and \cite{gr15}. By expressing the Green function of $Z$ through function $h$ we obtain estimates of the probability that $Z$ does not die upon exiting the interval $(0,R)$, as well as estimates of the expected exit time of $Z$ from the same interval in terms of $h$. These results can be also found in a recent paper \cite{gr15}, where they have been considered in a similar setting. 

Using these results, as well as sharp two-sided Green function estimates for $Z^{(a,b)}$ obtained through the conditional gauge theorem, we arrive to the main results of this paper by applying standard methods from \cite{onedimBM} and \cite{vondra}.

The paper is composed as follows. In Section \ref{prelim} we recall some basic results for a one-dimensional subordinate Brownian motion and consider several properties of function $h$. Applying these results, in Section \ref{properties} we prove several properties of the first exit time of $Z$ from a finite interval $(0,R)$. In Section \ref{green} we prove that process $Z$ killed outside of a finite interval $(a,b)$, $0<a<b$, can be obtained from the killed censored process $Y^{(a,b)}$ by a combination of a discontinuous and continuous Feynman-Kac transform and show that the Green functions for $X^{(a,b)}$ and $Z^{(a,b)}$ are comparable. Finally, in Section \ref{BHP} we give the proof of the Harnack inequality and boundary Harnack principle for $Z^{(a,b)}$.

\section{Preliminaries}\label{prelim}
Let $\phi:(0,\infty)\to(0,\infty)$ be a complete Bernstein function, $\phi\in\mathcal{CBF}$, with killing term and drift zero, i.e.~
\[
\phi(\lambda)=\int_0^\infty (1-e^{-\lambda t})\nu(t)dt,
\]
where the L\' evy measure satisfies the condition $\int_0^\infty (1\wedge t)\nu(t)dt<\infty$. Note that every Bernstein function $\phi$ satisfies the condition
\begin{equation}
  1\wedge\lambda\leq\frac{\phi(\lambda r)}{\phi(r)} \leq 1 \vee\lambda,\quad \lambda, r>0.\label{eq:phi}
\end{equation}
Let $W=(W_t)_{t\geq 0}$ be a $1$-dimensional Brownian motion and $S=(S_t)_{t\geq 0}$ a subordinator independent of $W$ with the Laplace exponent $\phi$, that is
\begin{equation*}
\mathbb E[e^{-\lambda S_t}]=e^{-t\phi(\lambda)},\quad t\geq 0,\,\lambda>0.
\end{equation*}
Define a one-dimensional subordinate Brownian motion $X=(X_t)_{t\geq 0}$ by $X_t=W_{S_t}$. Then $X$ is a L\' evy process with the characteristic exponent 
\[
\psi(x)=\phi(x^2),\quad x\in\mathbb R
\]   
and a decreasing L\' evy density
\[
j(x)=\int_0^\infty (4\pi s)^{-1/2}e^\frac{-x^2}{4s}\nu(s)ds,\, x\in\mathbb R.
\]

Furthermore, we will only consider the case when 0 is regular for itself, i.e. when
\[
\mathbb P_0(\sigma_{\{0\}}=0)=1,
\]
where $\sigma_{B}=\inf\{t>0:X_t\in B\}$ is the first hitting time of $B\in\B(\R)$ for $X$. By \cite[Lemma 3.1]{yano} this holds if and only if the Kesten-Bretagnolle condition is satisfied, that is if
\begin{equation}\label{eq:kesten}
\int_1^\infty\frac{1}{\phi(\lambda^2)}d\lambda<\infty.
\end{equation}
By \cite[Theorem II.16]{bertoin} there exists a bounded and continuous density $u^q$ of the $q$-resolvent
\[
U^qf(x)=\int_0^\infty e^{-qt}\EE_x[f(X_t)]dt=\int_{\R}f(x)u^q(x)dx
\]
of the form
\[
u^q(x)=\int_0^\infty e^{-qt}p_t(x)dt=\int_0^\infty e^{-qt}\frac{1}{2\pi}\int_{\R}e^{-i\lambda x}e^{-t\phi(\lambda^2)}d\lambda ~dt=\frac1{2\pi}\int_\R\frac{\cos(\lambda x)}{q+\phi(\lambda^2)}d\lambda.
\]
Since the transition density $p_t(x)$ is decreasing in $x$ it follows that $u^q$ is decreasing as well. 

\begin{Def} A Borel measurable function $f$ on $\mathbb R$ is harmonic on a Borel set $D$ with respect to a Markov process $X$  if for every bounded open set $B\subset \overline{B}\subset D$
\begin{equation}\label{eq:harmonic}
 f(x)=\mathbb E_x[f(X_{\tau_B})],\quad x\in B
\end{equation}
where $\tau_B=\inf\{t>0:X_t\not\in B\}$ is the first exit time of $X$ from $B$. If \eqref{eq:harmonic} holds also for $D$ in place of $B$, we say that $f$ is regular harmonic on $D$.
\end{Def}

Here we assume that the expectation in \eqref{eq:harmonic} is finite,  $X_{\infty}=\partial$, where $\partial$ is the so-called cemetery point and that $f(\partial)=0$. 

Define $h:\mathbb R\to[0,\infty)$ as 
\begin{equation}\label{killed:h}
\begin{aligned}
h(x)&=\lim\limits_{q\downarrow0}(u^q(0)-u^q(x))=\frac 1{\pi}\int_0^{\infty} \frac{1-\cos(\lambda x)}{\phi(\lambda^2)}d\lambda.
\end{aligned}
\end{equation}
The function $h$ is symmetric and since $u^q$ is decreasing, $h$ is increasing on $[0,\infty)$. 

Let $X^0$ be the process $X$ killed at 0 and $Z=(Z_t)_{t\geq 0}$ the absolute value of that process,
\[
 Z_t=
\left\{\begin{array}{c l}
 |X_t|,& t<\sigma_{\{0\}}\\
\partial, & t\geq\sigma_{\{0\}}
 \end{array},\quad t\geq 0.\right.
\] 
Since 0 is not polar, that is $\mathbb P_x(\sigma_0<\infty)>0$ for all $x\in\mathbb R$, $X^0$ is a proper subprocess of $X$. Also, if $X$ is recurrent then by \cite[Theorem 3.1]{yano} $\mathbb P_x(\sigma_0<\infty)=1$ for all $x\in\R$. By \cite[Theorem 1.1]{yano} $h$ is harmonic for the process $X^0$ on $\R\setminus\{0\}$ and since it is symmetric, it is also harmonic for $Z$ on $(0,\infty)$. Let $G^{X^0}(x,dy)$ and $G^{Z}(x,dy)$ be the Green measures for $X^0$ and $Z$ respectively. Note that for every $x>0$ and $A\in\mathcal B(0,\infty)$ 
\begin{align*}
 G^Z(x,A)&=\int_0^\infty (\P_x(X_t^0\in A)+ \P_x(-X_t^0\in A))dt=\int_A\left(G^{X^0}(x,y)+G^{X^0}(x,-y)\right)dy
\end{align*}
and thus the Green function of $Z$ is equal to 
\begin{equation}
G^{Z}(x,y)=G^{X^0}(x,y)+G^{X^0}(x,-y).\label{eq:greenZ_1} 
\end{equation}
Furthermore, the Green function $G^{X^0}$ of $X^0$ can be represented in terms of the function $h$. By \cite[Corrolary II.18]{bertoin} it follows that $\EE_x\left[e^{-q\sigma_{\{0\}}}\right]=\frac{u^q(x)}{u^q(0)}$ and therefore
\begin{align*}
G^{X^0}(x,y)&=\lim\limits_{q\downarrow0}u_0^q(x,y)=\lim\limits_{q\downarrow0}\left(u^q(x,y)-\EE_x\left[e^{-q\sigma_{\{0\}}}\right]u^q(0,y)\right)\\
&\overset{\eqref{killed:h}}{=}-h(y-x)+h(x)+h(y)-\kappa h(x)h(y),
\end{align*}
where $\kappa=\left(\frac 1 \pi\int_0^\infty\frac 1{\phi(\lambda^2)}d\lambda\right)^{-1}\in [0,\infty)$. Note that by the Chung-Fuchs type criteria for recurrence and \eqref{eq:kesten} $\kappa=0$ if and only if $X$ is recurrent. It follows from \eqref{eq:greenZ_1} that
\begin{equation}
G^{Z}(x,y)=2h(x)+2h(y)-h(y-x)-h(y+x)-2\kappa h(x)h(y),\quad x,y>0.\label{eq:greenZ}
\end{equation}
 
The following lemma is implied by \cite[Proposition 2.2, Proposition 2.4]{gr15}, which establish similar bounds for $G^{X^0}$ in terms of $h$.

\begin{Lem}\label{lemmaGz}For every $x,y>0$, $G^Z(x,y)\leq 4h(x\wedge y).$ If $X$ is recurrent then $G^Z(x,y)\geq h(x\wedge y)$, for every $x,y>0$. 
\end{Lem}
\proof
First we show that $h$ is a subadditive function on $\mathbb R$. Since $h$ is symmetric it follows that 
\begin{align*}
h(x)+h(y)-h(x+y)&=h(-x)+h(y)-h(x+y)\geq G^{X^0}(-x,y)\geq 0.
\end{align*}
By \eqref{eq:greenZ} and subadditivity of $h$ for $0<x<y$ we get
\begin{align*}
G^{Z}(x,y)&\leq 2h(x)+2h(y)-h(y-x)-h(y+x)-2\kappa h(x)h(y)\\
&\leq 2h(x)+h(x)+h(y-x)+h(-x)+h(y+x)-h(y-x)-h(y+x)=4 h(x). 
\end{align*}
Since $h$ is increasing on $(0,\infty)$, when $\kappa=0$ it follows that
\begin{align*}
G^{Z}(x,y)\geq h(x)+h(y)-h(y-x)\geq h(x).
\end{align*}
\qed
Throughout this paper we will assume that the complete Bernstein function $\phi$ satisfies the following global scaling condition\\
{\bf (H):} There exist constants $a_1,a_2>0$ and $\frac 1 2<\delta_1\leq\delta_2<1$ such that
\begin{equation*}
 	a_1\lambda^{\delta_1}\leq\frac{\phi(\lambda r)}{\phi(r)}\leq
 	a_2\lambda^{\delta_2}, \quad\lambda\geq 1, r>0. 
\end{equation*}

Note that since $\delta_1>\frac 1 2$ the regularity condition \eqref{eq:kesten} is satisfied and $\kappa=0$, i.e. $X$ is recurrent and 0 is regular for itself. 

We will use the following estimate of $h$ in terms of the characteristic function $\psi$ several times in the following chapter, see also \cite[Lemma 2.14]{gr15}.
\begin{Lem}\label{lemma1} There exists a constant $c_1>1$ such that for all $x>0$
$$c_1^{-1}\frac 1{x\psi\left(\frac 1 x\right)}\leq h(x)\leq c_1 \frac 1{x\psi\left(\frac 1 x\right)}.$$
\end{Lem}
\proof
For every $x\in\R$ it follows that
\begin{align*}
 h(x)&\leq\frac1\pi\int_0^\infty \left(\frac{\xi^2x^2}2\wedge 2\right)\frac 1{\psi(\xi)}d\xi= \frac{x^2}{2\pi}\int_0^{\frac 2 x}\frac {\xi^2}{\phi(\xi^2)}d\xi+\frac2\pi\int_{\frac 2 x}^\infty\frac 1{\phi(\xi^2)}d\xi\\
&\overset{\HH}{\leq} \frac{x^2}{2\pi}\frac{a_2}{\phi(4x^{-2})}\left(\frac 2 x\right)^{2\delta_2}\int_0^{\frac 2 x}\xi^{2-2\delta_2}d\xi+\frac{2x^{-2\delta_1}}{a_1\pi\phi\left(4 x^{-2}\right)}\int_{\frac 2 x}^\infty \xi^{-2\delta_1}d\xi\leq\tilde c_1\frac{1}{x\phi\left(x^{-2}\right)}.
\end{align*}
On the other hand, denoting the Fourier transform operator by $\mathcal F$ we have
\begin{align*}
 h(x)&=\frac1{2\pi}(\F\frac 1{\psi}(0)-\F\frac 1{\psi}(x))\geq\frac 1 {4\pi}\int_0^\infty\left(\frac{\xi^2x^2}4\wedge1\right)\frac 1{\psi(\xi)}d\xi= \frac{x^2}{\pi}\int_0^{\frac 2 x}\frac{\xi^2}{\phi(\xi^2)}d\xi+\frac1\pi\int_{\frac 2 x}^\infty\frac 1{\phi(\xi^2)}d\xi\\
 &\overset{\HH}{\geq} \frac{x^2}{4\pi}\frac{a_1}{\phi(4x^{-2})}\left(\frac 2 x\right)^{2\delta_1}\int_0^{\frac 2 x}\xi^{2-2\delta_1}d\xi+\frac{2x^{-2\delta_2}}{a_2\pi\phi\left(4 x^{-2}\right)}\int_{\frac 2 x}^\infty \xi^{-2\delta_2}d\xi\geq \tilde c_2\frac{1}{x\phi\left(x^{-2}\right)}.
\end{align*}
\qed

From the previous lemma and \HH it follows that $h$ also satisfies a global scaling condition, i.e.~there exist constants $d_1,d_2>0$ such that
\begin{equation}
 d_1\lambda^{2\delta_1-1}\leq\frac{h(\lambda t)}{h(t)}\leq d_2\lambda^{2\delta_2-1},\quad\forall \lambda\geq1,\,t>0. \label{eq:H3_h}
\end{equation}
 
This condition implies the following lower bound for the Green function $G^Z_{(0,R)}$ of the killed process $Z^{(0,R)}$, $R>0$, which is obtained similarly as in \cite[Lemma 4.2]{gr15}. We omit the proof.

\begin{Lem}\label{lemmaGz2}
There exist $\lambda_1\in\left(0,\frac 1 2\right)$ and $\lambda_2>0$ such that for every $R>0$
$$ G^Z_{(0,R)}(x,y)\geq \lambda_2h(R),\,\, x,y\in(0,\lambda_1 R).$$
\end{Lem}

\section{Properties of the exit time of $Z$ from a finite interval}\label{properties}

Denote by $\sigma_0:=\sigma_{\{0\}}$ the lifetime of $Z$ and $\tau_{(0,R)}=\inf\{t>0:Z_t\not\in(0,R)\}$ the first exit time of $Z$ from $(0,R)$, $R>0$. The following probability estimate that $Z$ does not die upon exiting $(0,R)$ was also obtained in \cite[Proposition 2.7]{gr15}.

\begin{Lem}\label{killed:lem1} For every $R>0$ and $x\in(0,R)$ 
 \begin{equation*}
\frac 1 8\frac{h(x)}{h(R)}\leq \P_x\left(\tau_{(0,R)}<\sigma_0\right)\leq \frac{h(x)}{h(R)}.
\end{equation*}
 \end{Lem}
\proof
First we prove the right inequality. For $\varepsilon>0$ by harmonicity of $h$ on $(0,\infty)$,
$$h(x)=\EE_x\left[h\left(Z_{\tau_{(\varepsilon,R)}}\right)\right]=\EE_x\left[h\left(Z_{\tau_{(\varepsilon,R)}}\right):\tau_{(\varepsilon,R)}<\sigma_{\{0\}}\right].$$
Since $h$ is continuous and $h(0)=0$ by the dominated convergence theorem and quasi-left continuity of $Z$ it follows that $h$ is regular harmonic for $Z$ on $(0,R)$,
\begin{align*}
h(x)&=\lim_{\varepsilon\to 0}\EE_x\left[h\left(Z_{\tau_{(\varepsilon,R)}}\right):\tau_{(\varepsilon,R)}<\sigma_{\{0\}}\right]=\EE_x\left[h\left(Z_{\tau_{(0,R)}}\right):\tau_{(0,R)}<\sigma_{\{0\}}\right]. 
\end{align*}
Since $h$ is increasing it follows that
\begin{align*}
h(x)&=\int_R^\infty h(y)\P_x\left(Z_{\tau_{(0,R)}}\in dy :\tau_{(0,R)}<\sigma_0\right)\geq h(R)\P_x\left(\tau_{(0,R)}<\sigma_0\right). 
\end{align*}

For the other inequality, by continuity and harmonicity of the Green function $G^Z(\cdot,2R)$ on $(\varepsilon,R)$ and Lemma \ref{lemmaGz}, it follows that
\begin{align*}
h(x)&\leq G^Z(x,2R)=\lim_{\varepsilon\to 0}\EE_x\left[G^Z(Z_{\tau_{(\varepsilon,R)}},2R)\right]=\int_R^\infty G^Z(z,2R)\P_x(Z_{\tau_{(0,R)}}\in dz)\\
&\leq 4h(2R) \P_x(\tau_{(0,R)}<\sigma_0)\overset{\eqref{eq:phi},\eqref{killed:h}}{\leq}  8h(R)\P_x(\tau_{(0,R)}<\sigma_0).
\end{align*}  
\qed

The following estimate for the tail distribution function of the lifetime of $Z$ was proven in \cite[Corollary 3.5.]{gr15}. Under additional assumptions it is also possible to obtain estimates of the derivatives of the tail distribution with respect to the time component. For more detail see \cite{kwasnicki}. 

\begin{Lem}\label{killed:lem2'}
 If there exist $a_1>0$ and $\delta_1\in(0,1]$ such that $\phi(\lambda t)\geq a_1\lambda^{\delta_1}\phi(t)$ hold for all $\lambda\geq1$ and $t>0$, then there exists a constant $c_2=c_2(n,\phi)$ such that   
 \begin{equation*}
  c_2^{-1}\frac{h(x)}{h\left(1/\psi^{-1}\left(\frac 1 t\right)\right)}\leq \P_x(\sigma_0>t)\leq c_2\frac{h(x)}{h\left(1/\psi^{-1}\left(\frac 1 t\right)\right)}
 \end{equation*}
for every $x\neq 0$ and $t>0$ such that $t\psi(\frac 1 x)\geq 1$. 
\end{Lem}

Using this result we can easily derive the following estimates for the expected exit time of $Z$ from interval $(0,R)$ in terms of the function $h$. 

\begin{Lem}\label{killed:lem3}
There exists a constant $c_3=c_3(R,\phi)>0$ such that
 \begin{align*}
  &\text{(i)}\quad\quad\EE_x\left[\tau_{(0,R)}\right]\leq 4Rh(x),\quad 0<x<R\\  
  &\text{(ii)}\quad\quad\EE_x\left[\tau_{(0,R)}\right]\geq c_3h(x),\quad \text{for }x\text{ small enough}.
 \end{align*}
\end{Lem}

\proof
(i) By Lemma \ref{lemmaGz} 
 \begin{align*}
  \EE_x\left[\tau_{(0,R)}\right]=\int_0^R G_{(0,R)}^Z(x,y)dy\leq \int_0^R 4h(x) dy=4Rh(x)
  \end{align*}
(ii) For the other inequality note that for all $t>0$
\begin{align*}
 \P_x(\sigma_0>t)&=\P_x(\sigma_0>t,\tau_{(0,R)}\geq\sigma_0)+\P_x(\sigma_0>t,\tau_{(0,R)}<\sigma_0)\\
 &\leq \P_x(\tau_{(0,R)}>t)+\P_x(\tau_{(0,R)}<\sigma_0)\leq \frac{\EE_x\left[\tau_{(0,R)}\right]}t+\P_x(\tau_{(0,R)}<\sigma_0),
\end{align*}
where the last term follows from Markov's inequality. Hence, by Lemma \ref{killed:lem1}, Lemma \ref{killed:lem2'} and Lemma\ref{lemma1}, if $t\psi\left(\frac 1 x\right)>1$ there exists a constant $\tilde c_1>0$ such that
\begin{align}
 \EE_x\left[\tau_{(0,R)}\right]&\geq t\left(\P_x(\sigma_0>t)-\P_x(\tau_{(0,R)}<\sigma_0)\right)\geq c_1t\frac{h(x)}{h\left(1/\psi^{-1}\left(\frac 1 t\right)\right)}-t\frac{h(x)}{h(R)}\nonumber\\
 &\geq\left(\frac{c_1\tilde c_1}{\psi^{-1}\left(\frac 1 t\right)}-\frac t{h(R)}\right)h(x)=f_{R}(t)h(x).
\end{align}
Note that by {\bf (H)} there exists a constants $\tilde c_2>0$ such that for all $t\leq 1$ 
$$f_{R}(t)\geq c_2\psi^{-1}(1) t^{\frac{-1}{2\delta_1}}-\frac t{h(R)},$$
so there exists $t_0=t_0(\phi,R)\in(0,1)$ such that $f_R(t)>0$ for all $t<t_0$. Therefore, 
\[
\EE_x[\tau_{(0,R)}]\geq f_R(t_0)h(x), \text{ for all }x<\frac1{\psi^{-1}(\frac{1}{t_0})}.
\]\qed

\section{Green function estimates for $Z^{(a,b)}$}\label{green}

Let $X^{(a,b)}$ and $Z^{(a,b)}$ be processes $X$ and $Z$ killed outside of interval $(a,b)$, $0<a<b$. In this section we obtain sharp bounds on the Green function $G^Z_{(a,b)}$ by comparing it to the Green function of $X^{(a,b)}$.

Let $Y$ be the process obtained from $X^{(0,\infty)}$ through the Feynman-Kac transform with respect to the positive continuous additive functional $A_\kappa$ with potential $\kappa_{(0,\infty)}(x)=\int_{-\infty}^0j(|x-y|)dy$, i.e.
\begin{equation}
\EE_x[f(Y_t)]=\EE_x\left[e^{A_\kappa(t)}f(X^{(0,\infty)}_t)\right]=\EE_x\left[e^{\int_0^t\kappa_{(0,\infty)}(X_s^{(0,\infty)})ds}f(X^{(0,\infty)}_t)\right]\label{eq:resurrect}
\end{equation}
for every bounded Borel function $f$ on $(0,\infty)$. We call $Y$ the resurrected (censored) process on $(0,\infty)$ corresponding to $X$, see \cite{cen} for a study of the censored process corresponding to a symmetric $\alpha$-stable L\' evy process, $\alpha\in(0,2)$. 

From the representation of Beurling-Deny and LeJan, the jumping measure associated with the Dirichlet form $(\E^Z,\F^Z)$ corresponding to the process $Z$ has a density equal to
\begin{equation*}
 i(x,y)=j(|x-y|)+j(|x+y|).
\end{equation*}
The Dirichlet forms corresponding to the processes $X^{(a,b)}$, $Y^{(a,b)}$ and $Z^{(a,b)}$ are therefore equal to
\begin{equation}\label{killed:dirichlet}
\begin{aligned}
 &\E^{X^{(a,b)}}(u,u)=\frac 1 2\int_a^b\int_a^b(u(x)-u(y))^2j(|x-y|)dy dx+\int_a^b u(x)^2\kappa_1(x)dx\\
 &\E^{Y^{(a,b)}}(u,u)=\frac 1 2 \int_a^b\int_a^b(u(x)-u(y))^2j(|x-y|)dy dx+\int_a^b u(x)^2\kappa_2(x)dx\\
 &\E^{Z^{(a,b)}}(u,u)= \frac 1 2\int_a^b\int_a^b(u(x)-u(y))^2i(x,y)dy dx+\int_a^b u(x)^2\kappa_3(x)dx,
\end{aligned}
\end{equation}
where the killing densities $\kappa_1$, $\kappa_2$ and $\kappa_3$ are of the form
\begin{align*}
 \kappa_1(x)=\int_{(a,b)^c}j(|x-y|)dy,\quad 
\kappa_2(x)=\int_{(0,\infty)\setminus (a,b)}j(|x-y|)dy,\quad
 \kappa_3(x)=\int_{(0,\infty)\setminus (a,b)}i(x,y)dy.
\end{align*}
Note that $Y^{(a,b)}$ can be obtained from $X^{(a,b)}$ by creation through the Feynman-Kac transform at rate $\kappa_{(0,\infty)}$. Therefore, by \cite[Lemma 3.4]{gauge} we can relate the Green functions of processes $X^{(a,b)}$ and $Y^{(a,b)}$ through a conditional gauge function $u_1(x,y)=\EE_x^y\left[e^{A_{\kappa}(\zeta^X_{(a,b)})} \right]$ by
$$G_{(a,b)}^Y(x,y)=u_1(x,y)G_{(a,b)}^X(x,y).$$
Here $\zeta^X_{(a,b)}=\inf\{t>0:X_t\not\in(a,b)\}$ is the lifetime of $X^{(a,b)}$ and $\P_x^y$ denotes the probability measure of the $G^X_{(a,b)}(\cdot,y)$-conditioned process starting from $x$, i.e.~the process with transition probability
\begin{equation}\label{eq:conditioned}
 p^y_t(x,z)=\frac{G^X_{(a,b)}(z,y)}{G^X_{(a,b)}(x,y)}p_t^{X^{(a,b)}}(x,z).
\end{equation}
Next we recall the definition of the Kato class $S_{\infty}$ from \cite{gauge}.

\begin{Def} Let $X$ be a transient Hunt process with the Green function $G$. A nonnegative Borel function $\kappa$ is said to be of the Kato class $S_\infty(X)$ if for any $\varepsilon>0$ there is a Borel set $K$ of finite measure and a constant $\delta>0$ such that
\begin{equation*}
\sup_{x,z\in \R^n}\int_{K^c\cup B}\frac{G^X(x,y)G^X(y,z)}{G^X(x,z)}\kappa(y)dy<\varepsilon 
\end{equation*}
for all measurable sets $B\subset K$ such that $\lambda(B)<\delta$.
\end{Def}
By the conditional gauge theorem \cite[Theorem 3.3]{gauge} if $\kappa_{(0,\infty)}\in S_\infty(X^{(a,b)})$ the conditional gauge function $u_1$ is bounded between two positive numbers. The key ingredient in showing $\kappa_{(0,\infty)}\in S_\infty(X^{(a,b)})$ is the following Green function estimate for $G^X_{(a,b)}$ from \cite[Corollary 7.4 (ii)]{plms14}. Let $\Phi(x):=\frac 1{\phi(x^{-2})}$, $\delta(x):=\text{dist}(x,(a,b)^c)$ and $a(x,y):=\Phi(\delta(x))^\frac 1 2\Phi(\delta(y))^\frac 1 2$. There exists a constant $c_4>1$ such that for every $x,y\in (a,b)$  
\begin{equation}
c_4^{-1}\left(\frac{a(x,y)}{\Phi^{-1}(a(x,y))}\wedge\frac{a(x,y)}{|x-y|}\right)\leq G_{(a,b)}^X(x,y)\leq c_4\left(\frac{a(x,y)}{\Phi^{-1}(a(x,y))}\wedge\frac{a(x,y)}{|x-y|}\right).\label{killed:G}
\end{equation}
Also, {\bf (H)} implies that $\Phi^{-1}$ satisfies the following scaling condition: there exists a constant $c_5>1$ such that for all $0<r\leq R<\infty$
\begin{equation}
c_5^{-1}\left(\frac r R\right)^{1/(2\delta_1)}\leq\frac{\Phi^{-1}(r)}{\Phi^{-1}(R)}\leq c_5\left(\frac r R\right)^{1/(2\delta_2)}.\label{eq:H4}
\end{equation}

\begin{Thm}\label{gaugeXY}
Let $X$ be a recurrent subordinate Brownian motion with Laplace exponent of the subordinator $\phi\in\mathcal{CBF}$ satisfying {\bf (H)} with $\delta_1>\frac 1 2$. Then the function $\kappa_{(0,\infty)}$ is in Kato class $S_\infty(X^{(a,b)})$. Therefore, the pair $(X^{(a,b)},\kappa_{(0,\infty)})$ is conditionally gaugeable and consequently the Green functions $G^X_{(a,b)}$ and $G^Y_{(a,b)}$ are comparable. 
\end{Thm}
\proof
Let $\varepsilon>0$. From \eqref{killed:G} we get the following 3G inequality,
\begin{align}
 A(x,y,z):= \frac{G_{(a,b)}^X(x,y)G_{(a,b)}^X(y,z)}{G_{(a,b)}^X(x,z)}\leq \frac{\tilde c_1 \Phi(\delta(y))(|x-z|\vee \Phi^{-1}(a(x,z)))}{(|x-y|\vee \Phi^{-1}(a(x,y)))(|y-z|\vee \Phi^{-1}(a(y,z)))}\label{killed:3G}
\end{align}
for some $\tilde c_1>0$. First note that by \eqref{eq:phi} and \eqref{eq:H4} if $\delta(y)\leq 2\delta(x)$ then
\begin{align*}
 \Phi^{-1}(a(x,y))&\geq \Phi^{-1}\left(\frac 1 4\Phi\left(\delta(y)\right)^{1/2}\Phi(\delta(y))^{1/2}\right)\geq c_5^{-1}2^{-\frac{1}{\delta_1}}\delta(y).
\end{align*}
Since $\delta(y)\leq 2(\delta(x)\vee |x-y|)$ it follows that
$|x-y|\vee \Phi^{-1}(a(x,y))\geq \left( \frac 1 2\wedge c_5^{-1}2^{-1/\delta_1}\right)\delta(y).$
Combining this inequality with \eqref{killed:3G} we arrive to
$A(x,y,z)\leq  \tilde c_1 \left( 4\vee c_5^{2}2^{2/\delta_1}\right)(b-a)\frac{\Phi(\delta(y))}{\delta(y)^2}$. Next, let $\tilde c_2=\tilde c_1 \left( 4\vee c_5^{2}2^{2/\delta_1}\right)(b-a)$ and $A=[a,a+\eta]\cup[b-\eta,b]$ for some $0<\eta<(b-a)\wedge 1$. It follows that
\begin{align*}
 \sup_{x,z\in(a,b)}\int_{A}A(x,y,z)dy&\leq 2\tilde c_2\int_0^\eta\frac{\Phi(s)}{s^2}ds\overset{{\bf (H)}}{\leq}\frac{ 2\tilde c_2}{a_1\phi(1)}\int_0^{\eta}\frac{s^{2\delta_1}}{s^2}ds= \tilde c_3\eta^{2\delta_1-1}.
\end{align*}
Therefore, by choosing $\eta$ small enough and $K:=[a+\eta,b-\eta]$ we get to
$ \displaystyle{\sup_{x,z\in(a,b)}\int_{K^c}A(x,y,z)dy< \frac\varepsilon 2}.$
The function $s\mapsto\frac{\Phi(s)}{s^2}$ is continuous on $[\eta,\frac{b-a}2]$ so there exists a constant $M>0$ such that $ \displaystyle{\sup_{x,z\in(a,b)}\int_{B}A(x,y,z)dy< \frac\varepsilon 2}$ for all $B\subset K$ such that $\lambda(B)<\delta:=\frac{\varepsilon}{2\tilde c_2M}$. Since $\kappa_{(0,\infty)}$ is bounded on $(a,b)$, this is enough to conclude that $\kappa_{(0,\infty)}\in S_\infty(X^{(a,b)})$.
\qed

Next, we associate the Green functions for processes $Y^{(a,b)}$ and $Z^{(a,b)}$. Since 
$$\E^{Z^{(a,b)}}(u,u)=\E^{Y^{(a,b)}}(u,u)+\int_a^b\int_a^b(u(x)-u(y))^2F(x,y)j(|x-y|)dy dx+\int_a^b u(x)^2q(x)dx,$$
where $F(x,y)=\frac{j(|x+y|)}{j(|x-y|)}$ and $q=\kappa_3-\kappa_2$, $Z^{(a,b)}$ can be obtained from $Y^{(a,b)}$ through the Feynman-Kac transform driven by a discontinuous additive functional 
\begin{equation}
A_{q+F}(t)= \int_0^t q(Y^{(a,b)}_s)ds+\sum_{s\leq t} F(Y^{(a,b)}_{s-},Y^{(a,b)}_s).\label{killed:AqF}
\end{equation}
By \cite[Lemma 3.9]{gauge} the ratio of Green functions $G_{(a,b)}^Z(x,y)$ and $G_{(a,b)}^Y(x,y)$ is equal to the gauge function $u_2(x,y)=\EE_x^y\left[ e^{A_{q+F}(\zeta^Y_{(a,b)})}\right]$ and $\zeta^Y_{(a,b)}=\inf\{t>0:Y_t\not\in(a,b)\}$ is the lifetime of $Y^{(a,b)}$ and $\mathbb P_x^y$ is the probability measure of the $G_{(a,b)}^Y(\cdot,y)$-conditioned process starting from $x$, see \eqref{eq:conditioned}. We recall the definition of the Kato class $A_\infty$. 

\begin{Def} Let $X$ be a transient Hunt process with values in $E\in \B(\R)$ with Green function $G$ and L\' evy system $(J,H)$, where $H_s\equiv s$. A bounded nonnegative function $F$ on $E\times E$ vanishing on the diagonal is said to be in the Kato class $A_{\infty}(X)$ if for any $\varepsilon > 0$ there is a Borel subset $K$ of finite measure and a constant $\delta> 0$ such that for every set $A=(K\times K)^c\cup (B\times E)\cup(E\times B)$ 
\begin{align*}
&\sup_{x,w\in E}\int_{A}\frac{G(x,y)G(z,w)}{G(x,w)}F(y,z)J(x,dy)dz<\varepsilon,
\end{align*}
where $B\subset K$ is a measurable set such that 
$\int_B\left(\int_{E}F(x,y)J(x,dy)\right)dx<\delta$.
\end{Def}

By \cite[Theorem 3.8]{gauge} the conditional gauge function $u_2$ is bounded between two positive constants when $q\in S_\infty(Y^{(a,b)})$ and $F\in A_{\infty}(Y^{(a,b)})$. This is shown by using \eqref{killed:G} similarly as in Theorem \ref{gaugeXY}, so we omit the proof.

\begin{Thm}\label{gaugeYZ} Let the assumptions from Theorem \ref{gaugeXY} hold and $A_{q+F}$ be the discontinuous additive functional for $Y^{(a,b)}$ from \eqref{killed:AqF}. Then $q\in S_\infty(Y^{(a,b)})$ and $F\in A_\infty(Y^{(a,b)})$ and consequently the Green functions of the processes $Y^{(a,b)}$ and $Z^{(a,b)}$ are comparable.
\end{Thm}

\section{Boundary Harnack principle for $Z$}\label{BHP}

The exit distribution of $Z_{\tau_{(a,b)}}$ starting from $x$ is equal to
\[
\P_x\left(Z_{\tau_{(a,b)}}\in B\right)=\int_B K^Z_{(a,b)}(x,z)dz,\quad x\in (a,b),\,B\in\B((0,\infty)\setminus[a,b]),
\]
where $K^Z_{(a,b)}$ is the Poisson kernel of $Z^{(a,b)}$ given by
\[
K^Z_{(a,b)}(x,z)=\int_a^b G^Z_{(a,b)}(x,y)i(y,z)dy,\quad x\in(a,b),\,z\in(0,\infty)\setminus[a,b]. 
\]
Recall that the process $Z$ can exit the interval $(a,b)$ only by jumping out, since by \cite[Theorem 1]{sztonyk}
\begin{align*}
\P_x\left(X_{\tau_{(a_1,a_2)}}=a_i\right)=0, \quad\, i=1,2
\end{align*}
for all $x\in(a_1,a_2)\subset\R$. Using the results from the previous sections we can similarly as in \cite[Section 4]{onedimBM} prove the Harnack inequality and boundary Harnack principle for nonnegative harmonic functions of process $Z^{(a,b)}$.

\begin{Thm}{\bf Harnack inequality}\label{killed:harnack}\\
Let $R > 0$ and $a\in(0,1)$. There exists a constant $c_6 = c_6(R,a,\phi) > 0$ such that for all $r \in (0,R)$ and every nonnegative function $u$ on $\R$ which is harmonic with respect to $Z$ in $(0, 3r)$,
$$u(x) \leq c_6 u(y), \text{  for all } x, y \in (ar, (3-a)r).$$
\end{Thm}

\proof
Let $b_1 = ar/2$, $b_2 = ar$, $b_3 = (3-a)r$ and $b_4 = (3-a/2)r$. By Theorem \ref{gaugeXY}, Theorem \ref{gaugeYZ} and \eqref{killed:G} the exists a $\tilde c_1=\tilde c_1(\phi, R)>1$ such that
\begin{align*}
\tilde c_1^{-1}\frac{a(x_i,y)}{\Phi^{-1}(a(x_i,y))\vee|x-y|}\leq  G^Z_{(b_1,b_4)}(x_i, y) \leq \tilde c_1\frac{a(x_i,y)}{\Phi^{-1}(a(x_i,y))\vee|x-y|},\quad i=1,2,
\end{align*}
for all $x_1, x_2 \in (b_2, b_3)$ and $y \in (b_1, b_4).$ Furthermore, note that
\begin{align*}
\frac {ar}2\leq \delta(x_i)\leq\frac{(3-a)r}2 \quad\text{ and }\quad\delta(y)\leq \frac {ar}4\, \Rightarrow\, |x_i-y|\geq \frac{ar}{4}.
\end{align*}
Therefore $\Phi^{-1}(a(x_i,y))\vee|x_i-y|$ is comparable to $r$, so by {\bf (H)} and \eqref{eq:H4} there exists a constant $\tilde c_2=\tilde c_2(R,a,\phi)>0$ such that
\begin{align*}
G^Z_{(b_1,b_4)}(x_1, y) &\leq \tilde c_2 G^Z_{(b_1,b_4)}(x_1, y)
\end{align*}
for all $x_1, x_2 \in (b_2, b_3)$ and $y \in (b_1, b_4).$ Consequently, we have
\begin{align*}
K^Z_{(b_1,b_4)}(x_1, z)\leq\tilde c_2 K^Z_{(b_1,b_4)}(x_2, z)
\end{align*}
for all $x_1, x_2 \in (b_2, b_3)$, $z \in [b_1, b_4]^c$. It follows that for all $x_1, x_2\in(ar, (3-a)r)$
\begin{align*}
u(x_1)=\EE_{x_1}\left[u\left(X_{\tau_{(b_1,b_4)}}\right)\right]=\int_{(b_1,b_4)^c} u(z)K^Z_{(b_1,b_4)}(x_1,z)\leq \tilde c_2\int_{(b_1,b_4)^c} u(z)K^Z_{(b_1,b_4)}(x_2,z)=\tilde c_2 u(x_2).
\end{align*}

\qed

\begin{Thm}{\bf Boundary Harnack principle}\\
Let $R > 0$. There exists a constant $c_7 = c_7(R,\phi) >
0$ such that for all $r \in (0,R)$, and every nonnegative function $u$ which is harmonic for $Z$ in $(0, 3r)$ and continuously vanishes at 0 it holds that
$$\frac{u(x)}{u(y)}\leq c_7\frac {h(x)}{h(y)},\quad x, y \in (0, \lambda_1 r),$$
where $\lambda_1$ is the constant from Lemma \ref{lemmaGz2}.
\end{Thm}
\proof
Let $x \in (0, \lambda_1r)$. Since $u$ is harmonic in $(0, 3r)$ and vanishes continuously at 0 we have
\begin{align*}
u(x) = &\lim_{\varepsilon\to0} \EE_x\left[u\left(Z_{\tau_{(\varepsilon,r)}}\right)\right]=\EE_x\left[u\left(Z_{\tau_{(0,r)}}\right)\right]=\EE_x\left[u\left(Z_{\tau_{(0,r)}}\right):Z_{\tau_{(0,r)}}\in(r,2r)\right]\\
&+\EE_x\left[u\left(Z_{\tau_{(0,r)}}\right):Z_{\tau_{(0,r)}}\geq 2r\right]=: u_1(x) + u_2(x).
\end{align*}
First note that
$$\frac{u(x)}{u(\lambda_1 r)}\leq\frac{u_1(x)}{u(\lambda_1 r)}+\frac{u_2(x)}{u_2(\lambda_1 r)}$$
and we estimate each term separately. By the previous Harnack inequality for $a=\frac{\lambda_1}2$ and Lemma \ref{killed:lem1} there exists a constant $\tilde c_1=\tilde c_1(R,\phi)>0$ such that 
\begin{align*}
u_1(x)&\leq \tilde c_1\EE_x\left[u\left(\lambda_1 r\right):Z_{\tau_{(0,r)}}\in(r,2r)\right]\leq \tilde c_1u(\lambda_1r)\P_x\left(\tau_{(0,r)}<\tau\right)\leq \tilde c_1u(\lambda_1r)\frac{h(x)}{h(\lambda_1r)}. 
\end{align*}
For the second term, since the L\' evy density $j$ of $X$ is decreasing it follows that
\begin{align*}
u_2(x)&=\int_0^r\int_{2r}^\infty u(z)G^Z_{(0,r)}(x,y)i(y,z)dz dy\leq \int_0^r G^Z_{(0,r)}(x,y)dy\int_{2r}^\infty u(z)(j(z-r)+j(z))dz\\
&= \EE_x[\tau_{(0,r)}]\int_{2r}^\infty u(z)(j(z-r)+j(z))dz\leq 4r h(x) \int_{2r}^\infty u(z)(j(z-r)+j(z))dz
\end{align*}
where the last line follows from Lemma \ref{killed:lem3}. By \cite[Theorem 3.4]{gubhp} there exists a constant $\tilde c_2=\tilde c_2(\phi)>0$ such that
\[
\tilde c_2^{-1} \frac{\phi(z^{-2})}{z}\leq j(z)\leq\tilde c_2 \frac{\phi(z^{-2})}{z},\quad z>0,
\] 
so by $\eqref{eq:phi}$ it follows that $j(z-r)\leq \tilde c_2 2^3 j(z)$ when $z\geq2r$ and therefore
\begin{align*}
u_2(x)\leq \tilde  4c_2^2(2^3+1) rh(x) \int_{2r}^\infty u(z)j(z)dz.
\end{align*}
On the other hand, by Lemma \ref{lemmaGz2}
\begin{align*}
u_2(x)&\geq\int_0^{\lambda_1r} G^Z_{(0,r)}(x,y)dy\int_{2r}^\infty u(z)(j(z)+j(z+r))dz\geq \lambda_2 h(\lambda_1 r)\lambda_1 r\int_{2r}^\infty u(z)j(z)dz.
\end{align*}
Therefore, it follows that 
\begin{equation}\label{killed:one}
\frac{u(x)}{u(\lambda_1 r)}\leq \left(\tilde c_1+\frac{4\tilde  c_2(1+2^3)}{\lambda_1\lambda_2}\right)\frac{h(x)}{h(\lambda_1 r)}.
\end{equation}
On the other hand, from \cite[Lemma 5.1]{vondra} for $p=\frac 13$ it follows that there exists a constant $\tilde c_3=\tilde c_3(\phi,R)>0$ such that for all $x\in(0,r)$ and $y\in(2r,3r)$
\begin{align*}
\int_{2r}^y K_{(0,s)}^Z(x,y)ds&\leq\int_{3r(1+1/3)/2}^y (K_{(-s,s)}^X(x,y)+K_{(-s,s)}^X(x,-y))ds\leq \frac{3\tilde c_3r}{\phi((3r)^{-2})}j(y)\overset{\eqref{eq:phi}}{\leq}  \frac{27\tilde c_3r}{\phi(r^{-2})}j(y).
\end{align*}

Now by applying \cite[Lemma 5.2 and Lemma 5.3]{vondra} for $U=B(0,2r)$ and $p=\frac 1 3$ it follows that
\[
 u(x)\leq \frac{\tilde c_4}{\phi(r^{-2})}\int_{2r}^\infty u(y) j(y)dy
\]
for some constant $\tilde c_4=\tilde c_4(\phi)>0$ and all $x\in(0,r)$. Furthermore by Lemma \ref{lemmaGz}
\begin{align*}
u_2(x)\geq\int_0^{\lambda_1r} G^Z_{(0,r)}(x,y)dy\int_{2r}^\infty u(z)(j(z)+j(z+r))dz\geq \lambda_2 h(x)\lambda_1 r\int_{2r}^\infty u(z)j(z)dz.
\end{align*}
By the last two displays, \eqref{eq:phi} and Lemma \ref{lemma1} we get the required inequality, i.e.
\begin{equation}\label{killed:two}
\frac{u(x)}{u(\lambda_1 r)}\geq\frac{u_2(x)}{u(\lambda_1 r)}\geq \frac{\lambda_1\lambda_2rh(x)}{\frac{\tilde c_4}{\phi(r^{-2})}}\geq \tilde c_5\frac{h(x)}{h(\lambda_1r)}.
\end{equation}
Combining \eqref{killed:one} and \eqref{killed:two} we get the statement of the theorem. 
\qed


\vspace{0.5cm}
Contact: \href{mailto:wagner@math.hr}{wagner@math.hr} 

\end{document}